\magnification=1200
\hsize=12.5cm
\vsize=19.5cm
\parskip 3pt plus 1pt minus 1pt
\parindent=6.4mm

\let\sl\it

\font\eightrm=cmr8
\font\eightbf=cmbx8
\font\eightsy=cmsy8
\font\eighti=cmmi8
\font\eightit=cmti8
\font\eighttt=cmtt8
\font\eightsl=cmsl8

\font\sixrm=cmr6
\font\sixbf=cmbx6
\font\sixsy=cmsy6
\font\sixi=cmmi6

\font\fourteenbf=cmbx10 at 14.4pt
\font\twelvebf=cmbx10 at 12pt

\font\tenmsa=msam10
\font\sevenmsa=msam7
\font\fivemsa=msam5
\newfam\msafam
  \textfont\msafam=\tenmsa 
  \scriptfont\msafam=\sevenmsa
  \scriptscriptfont\msafam=\fivemsa

\font\tenmsb=msbm10
\font\eightmsb=msbm8
\font\sevenmsb=msbm7
\font\fivemsb=msbm5
\newfam\msbfam
  \textfont\msbfam=\tenmsb
  \scriptfont\msbfam=\sevenmsb
  \scriptscriptfont\msbfam=\fivemsb
\def\Bbb{\fam\msbfam\tenmsb}

\catcode`\@=11
\def\eightpoint{%
  \textfont0=\eightrm \scriptfont0=\sixrm \scriptscriptfont0=\fiverm
  \def\rm{\fam\z@\eightrm}%
  \textfont1=\eighti \scriptfont1=\sixi \scriptscriptfont1=\fivei
  \def\oldstyle{\fam\@ne\eighti}%
  \textfont2=\eightsy \scriptfont2=\sixsy \scriptscriptfont2=\fivesy
  \textfont\itfam=\eightit
  \def\it{\fam\itfam\eightit}%
  \textfont\slfam=\eightsl
  \def\sl{\fam\slfam\eightsl}%
  \textfont\bffam=\eightbf \scriptfont\bffam=\sixbf
  \scriptscriptfont\bffam=\fivebf
  \def\bf{\fam\bffam\eightbf}%
  \textfont\ttfam=\eighttt
  \def\tt{\fam\ttfam\eighttt}%
  \def\Bbb{\fam\msbfam\eightmsb}%
  \textfont\msbfam=\eightmsb
  \def\Cal{\fam\Calfam\eightCal}%
  \textfont\Calfam=\eightCal
  \abovedisplayskip=9pt plus 2pt minus 6pt
  \abovedisplayshortskip=0pt plus 2pt
  \belowdisplayskip=9pt plus 2pt minus 6pt
  \belowdisplayshortskip=5pt plus 2pt minus 3pt
  \smallskipamount=2pt plus 1pt minus 1pt
  \medskipamount=4pt plus 2pt minus 1pt
  \bigskipamount=9pt plus 3pt minus 3pt
  \normalbaselineskip=9pt
  \setbox\strutbox=\hbox{\vrule height7pt depth2pt width0pt}%
  \let\bigf@ntpc=\eightrm \let\smallf@ntpc=\sixrm
  \normalbaselines\rm}
\catcode`\@=12

\def\bC{{\Bbb C}}

\def\bQ{{\Bbb Q}}
\def\bR{{\Bbb R}}

\font\tenCal=eusm10
\font\eightCal=eusm8
\font\sevenCal=eusm7
\font\fiveCal=eusm5
\newfam\Calfam
  \textfont\Calfam=\tenCal
  \scriptfont\Calfam=\sevenCal
  \scriptscriptfont\Calfam=\fiveCal
\def\Cal{\fam\Calfam\tenCal}

\def\cC{{\Cal C}}

\def\vep{\varepsilon}
\def\ddbar{\partial\overline\partial}
\def\\{\hfil\break}
\let\wt\widetilde
\let\wh\widehat

\def\hexnbr#1{\ifnum#1<10 \number#1\else
 \ifnum#1=10 A\else\ifnum#1=11 B\else\ifnum#1=12 C\else
 \ifnum#1=13 D\else\ifnum#1=14 E\else\ifnum#1=15 F\fi\fi\fi\fi\fi\fi\fi}
\def\msatype{\hexnbr\msafam}
\def\msbtype{\hexnbr\msbfam}
\mathchardef\compact="3\msatype62
\mathchardef\smallsetminus="2\msbtype72   
\mathchardef\subsetneq="3\msbtype28
\def\buildo#1\over#2{\mathrel{\mathop{\null#2}\limits^{#1}}}
\def\buildu#1\under#2{\mathrel{\mathop{\null#2}\limits_{#1}}}

\def\square{{\hfill \hbox{
\vrule height 1.453ex  width 0.093ex  depth 0ex
\vrule height 1.5ex  width 1.3ex  depth -1.407ex\kern-0.1ex
\vrule height 1.453ex  width 0.093ex  depth 0ex\kern-1.35ex
\vrule height 0.093ex  width 1.3ex  depth 0ex}}}

\let\Item\item
\def\item#1{\Item{$\rlap{\hbox{#1}}\kern\parindent\kern-5pt$}}

\def\today{\ifcase\month\or
January\or February\or March\or April\or May\or June\or July\or August\or
September\or October\or November\or December\fi \space\number\day,
\number\year}


\null
\vskip 2cm
\centerline{\fourteenbf Regularity properties of the}
\bigskip
\centerline{\fourteenbf degenerate Monge--Amp\`ere equations on}
\bigskip
\centerline{\fourteenbf compact K\"ahler manifolds}
\null
\vskip 10pt
\centerline{ Mihai P\u AUN }
\vskip 5pt
\centerline{ Institut  \'Elie Cartan, Nancy }
\medskip

\vskip 30pt

\noindent{\twelvebf 0. Introduction} 
\medskip

\noindent 
Let $(X, \omega)$ be a $n$--dimensional compact K\"ahler manifold, and 
let $\alpha\in \cC^\infty_{1,1}(X, \bR)$ be a closed 
$(1,1)$ form on $X$, such that:

\item {(1)} $\alpha\geq 0$ pointwise on $X$ and 
$\int_X\alpha^n > 0$;
\smallskip
\item {(2)} $\displaystyle {{\omega^n}\over {\alpha^n}}\in
L^{\varepsilon_0}(X, \omega)$ for some $\vep_0>0$.

\noindent A large class of such $(1,1)$--forms can be obtained as
follows: let $\pi:X\to Y$ be a generically finite map, and let $\omega_Y$ be a 
K\"ahler metric on $Y$; then
$\alpha:= \pi^*\omega_Y$ verify the conditions above.

\vskip 5pt

We recall that according to [3], a function
$\phi:X\to [-\infty, \infty)$ is called quasi-plurisubharmonic
(quasi-psh for short) if it is locally equal to the sum of a smooth 
function and a plurisubharmonic (psh) function.
Then there exist a constant
$C\in \bR$ such that $\sqrt{-1}\ddbar \phi\geq -C\omega$
in the sense of currents on $X$. We say that a function 
$\psi $ has {\it logarithmic poles} if for each open set
$U\subset X$ there exist a family of holomorphic functions
$(f_j^U)$ such that $\psi\equiv \sum_j\vert f_j^U\vert ^2$ modulo
$\cC^\infty(U)$; it is an important class of
quasi-psh functions.   
\vskip 5pt
In this setting, the aim of our note is to 
prove the next result.

\medskip
\noindent{\bf Theorem.} {\sl Let $(X, \omega)$ be a compact K\"ahler
manifold, and let $\alpha$ be a smooth $(1,1)$ form on $X$, having the
properties (1) and (2) above. Consider the quasi--psh functions
$\psi_1, \psi_2$, such that 
\item {(i)}
$\displaystyle 
\int_X\exp\bigl(p(\psi_1-\psi_2)\bigr)dV_\omega<\infty$ for some
$p> 1$;
\smallskip
\item {(ii)} $\int_X\alpha^n= \int_X\exp(\psi_1-\psi_2)dV_\omega$

\noindent Then for each $\gamma\in [0, 1)$, there exist $Y_\gamma\subset X$
such that the solution $\varphi$ of the equation

$$(\alpha+ \sqrt{-1}\ddbar \varphi)^n= 
\exp(\psi_1-\psi_2)\omega^n\leqno(*) $$
belongs to the H\"older class $\cC^{1, \gamma}(X\setminus Y_\gamma)$.
}
\medskip

\noindent The estimates for the norm of the solution above 
can be obtained from the proof, but since they are not very
enlighting, we have decided to skip them.

\vskip 5pt

\noindent Before stating the next result, we would like to 
recall the following conjecture, formulated by J.-P. Demailly--J. Koll\`ar
in [4].

\medskip

\noindent {\bf Conjecture {[4]}} {\sl Let $\psi$ be a psh function on the unit ball  
$B\subset \bC^n$, such that $\displaystyle \int_B\exp (-\psi)d\lambda< \infty$.
Then there exists a positive real number $\delta> 0$ such that 
$\displaystyle \int_{(\bC^n, 0)}\exp \bigl(-(1+ \delta)\psi\bigr)d\lambda< \infty$.
}

\medskip

\noindent A consequence of this conjecture 
would be the following statement.

\medskip
\noindent{\bf Conjecture.} {\sl Let $(X, \omega)$ be a compact K\"ahler
manifold, and let $\psi$be a quasi-psh function such that 
$\displaystyle \int_X\omega^n= \int_X\exp (-\psi)dV_\omega< \infty$. Then the
Monge--Amp\`ere equation
$$(\omega+ \sqrt{-1}\ddbar \varphi)^n= 
\exp(-\psi)\omega^n$$
has a unique continuous solution up to normalisation.}

\medskip
Now if the previous statement is correct, then by 
the Chern--Levine--Nirenberg inequalities
we get $\psi\exp(-\psi)\in L^1(X)$. 
Let us assume for simplicity that $\psi\leq -1$ on $X$; then
the function 
$\displaystyle \psi- \log (-\psi)$ is equally quasi-psh, so inductively
we would get $(-\psi)^k\exp (-\psi)\in L^1(X)$. Unfortunately, the 
$L^1$--integrability statements above seem to be 
more or less equivalent to the existence of the continuous solution in $(*)$.

\vskip 10pt

\noindent Anyway, an immediate consequence of the theorem above
in the next statement, which address a question proposed in the paper
[6].

\medskip
\noindent{\bf Corollary.} {\sl Let $(X, \omega)$ be a compact K\"ahler
manifold, and let $\alpha$ be a smooth $(1,1)$ form on $X$, with the
properties 1--3 above. Consider the quasi--psh functions
$\psi_1, \psi_2$, which are assumed to have logarithmic poles and satisfy the 
following conditions 
$\displaystyle \int_X\alpha^n=
\int_X\exp(\psi_1-\psi_2)dV_\omega<\infty$. 
Then there exist an analytic set $Y\subset X$
such that the solution $\varphi$ of the equation
$$(\alpha+ \sqrt{-1}\ddbar \varphi)^n=
\exp(\psi_1-\psi_2)\omega^n\leqno(4) $$
is smooth on $X\setminus Y$.
}

\medskip In the paper [6], the authors proved the corollary 
under the assumptions that the manifold $X$ is projective,
and the cohomology class $\{\alpha \}$ 
lies in the real Neron-Severi group of 
$X$ (see also the references in [6] for other results in this
direction).

\vskip 10pt

Let us outline the main steps in the proof of the theorem.
First of all, if the functions $(\psi_j)_{j=1,2}$ are regular enough, 
and if $\alpha$ is a genuine K\"ahler metric,
then the question was completely solved by S.-T. Yau
in [13]. Therefore, the natural idea under the hypothesis of the theorem is
to regularise the functions $(\psi_j)_{j=1,2}$, then use the result of Yau
in order to solve the equations with the regularised right hand side member,
and finally to take the limit. As it is well known, to carry out this programme 
we have to provide
uniform a-priori estimates for the solutions.

Now, if $\alpha$ is a K\"ahler metric, the $\cC^0$--estimates we need
were obtained by 
S. Kolodziej in [7] by using simple, tricky and elegant considerations in the
pluripotential theory. In a recent paper [6], Ph. Eyssidieux, V. Guedj and A. Zeriahi
proved that the methods of S. Kolodziej can be extended 
to cover the case where $\alpha$ is only 
semi-positive on $X$ and strictly positive at some
point of $X$; therefore,
the solution of the equation $(*)$ is known to be continuous.

In order to achieve further regularity, we would like to use 
the $\cC^2$--estimates in
the theory of Monge-Amp\`ere equations, but a new difficulty occur:
since the $(1,1)$--form $\alpha$ may be zero at some points of $X$, 
the quantity $tr_\alpha\ddbar f$ (clearly needed in the estimates)
will be unbounded, even for smooth functions $f$. 
This difficulty already appeared in H. Tsuji's paper [12]; there
he consider the case where
$X$ is projective, and $\alpha= c_1(L)$,
for some big line bundle $L$. 
He solved the problem by modifying the
$(1,1)$--form $\alpha$ within its cohomology class 
in order to get a strongly positive representative;
remark that this is possible by
the well-known fact (the Kodaira's lemma):
a big line bundle can be decomposed as a sum of effective and
ample $\bQ$--line bundles. Then he observed that the singularities 
which come into the picture via the effective part of the decomposition do not affect
in a significant way the usual $\cC^2$--estimates
(recently, the same circle of ideas were used in
[6] and [11]). 
 
In our case, the $(1,1)$--form
$\alpha$ do not correspond to a line bundle, and $X$ is not necessarily projective, 
but recall that a result of J.-P. Demailly and ourself [5]
show that if $\alpha$ is a semi-positive $(1,1)$--form such that
$\int_X\alpha^n> 0$, then there exist 
a function 
$\tau$ with at worst logarithmic poles, such that the current
$T:= \alpha+ \sqrt{-1}\ddbar \tau$ dominates a small multiple of the
K\"ahler metric $\omega$. Thus, we ``trade'' the smoothness of
$\alpha$ for the strong positivity of $T$, as in the case of line bundles; 
the only new phenomenon is that the poles of $\tau$ may not be of divisorial type.
 The existence of this current is crucial
for the regularity analysis, since it is the right substitute for
the Kodaira lemma quoted above.  

For the rest of the proof, we follow the classical approach
in the Monge--Amp\`ere theory
and we show that Y.-T. Siu's version of the second order estimates in [9]
can be adapted in our context to give the result.

\vskip 40pt

\noindent{\twelvebf 1. Regularization of currents and 
$\cC^0$--estimates} 
\medskip

\noindent 
As the title of this paragraph try to suggest,
we will collect here some facts about the regularization of
quasi-psh functions. We also recall some results concerning
the $\cC^0$--estimates for the Monge--Amp\`ere operators
which will be used later. The convention all over this  
paper is that we will use same letter ``$C$'' to denote
a generic constant, which may change from one line to another,
but it is independent of the pertinent parameters involved.
\vskip 5pt

Let $(X, \omega)$ be a $n$--dimensional compact K\"ahler 
manifold, and let $\psi$ be a quasi-psh function. By definition, 
there exist a constant $C> 0$ such that $C\omega+ \sqrt{-1}\ddbar \psi\geq 0$
on $X$. 
We recall the next result due to J.-P. Demailly,
on the regularization of $\psi$; in fact, the statement in [3] 
is much more precise, but
all we need is the following particular case.

\medskip
\noindent{\bf Theorem {\rm [3]}} {\sl There exists a family of smooth functions
$(\psi_\vep)\subset \cC^\infty(X)$ and a constant $C> 0$ such that:
\item {(i)} $\psi_\vep\to \psi$ in $L^1(X)$ as $\vep\to 0$, and $\psi_\vep \geq \psi- 1$ 
for all $0<\vep\ll 1$;
\smallskip
\item {(ii)} $C\omega+ \sqrt{-1}\ddbar \psi_\vep\geq 0$.}

\medskip
\noindent The functions $\psi_\vep$ are obtained by means of 
the flow of the Chern connection on the tangent bundle $T_X$; this can be seen as
a global version of the familiar local convolution by smoothing kernels.
We apply the previous regularization theorem to $\psi_1$ and $\psi_2$;
thanks to the fact that $\exp (\psi_1-\psi_2)$ is in $L^p$, for some $p> 1$,
by the above considerations we infer that $\exp (\psi_{1;\vep}-\psi_{2; \vep})\to
\exp (\psi_1-\psi_2)$ in $L^p$.
\vskip 5pt
We recall now the theorem of S.-T. Yau (see [13]), which will be used
(in direct or indirect manner) several times in this note.

\medskip
\noindent{\bf Theorem {\rm [13]}} {\sl Let $(X, \omega)$ be a compact K\"ahler manifold,
and let $dV$ be a smooth volume element, such that $\int_XdV= \int_X\omega^n$. Then there exist a 
smooth function $\varphi$, unique up to normalisation,
 such that $\omega + \sqrt{-1}\ddbar \varphi > 0$ and such that
$$(\omega + \sqrt{-1}\ddbar \varphi)^n= dV.$$}
\medskip

By hypothesis, $\alpha$ is a semi-positive 
$(1,1)$--form of positive top self-intersection, thus 
for each $\vep> 0$ we have  
$\alpha+ \vep \omega> 0$, and by the previous result,
there exist a function $\varphi_\vep\in \cC^\infty(X)$ such that 
$\alpha+ \vep\omega+ \sqrt{-1}\ddbar \varphi_\vep> 0$ on $X$ and which is a solution of the 
equation
$$(\alpha+ \vep\omega + \sqrt{-1}\ddbar \varphi_\vep)^n= (1+ \delta_\vep)
\exp (\psi_{1;\vep}-\psi_{2; \vep})\omega^n\leqno (5)$$
on $X$. We assume that $\displaystyle \int_X\varphi_\vep dV_\omega= 0$;
in the previous expression, the real numbers $\delta_\vep$
are normalisation constants, and we have $\delta_\vep\to 0$
as $\vep\to 0$.  

We want to prove that some subsequence of the family $(\varphi_\vep )$
converge to the solution of the equation (4), and that this limit has 
the regularity properties stated in the theorem. 
\vskip 5pt
A first step in this direction is provided by the next result, due to
S. Kolodziej, see [7].

\medskip
\noindent{\bf Theorem {\rm [7]}} {\sl Let $(X, \alpha)$ be a compact K\"ahler
manifold and let $(f_j)\subset \cC^\infty(X)$ be a sequence 
of functions on $X$, such that the following requirements are satisfied:
\item {(a)}
$\displaystyle \sup_j\Vert\exp(f_j)\Vert_{L^p(X)}< \infty$
for some $p> 1$; 
\smallskip
\item {(b)} $\int_X\exp(f_j)dV_\alpha= \int_X\alpha^n$ for $j\geq 1$. 

\noindent Then there exist a constant $C\in \bR$ such that
for each solution $\varphi_j$ of the equation
$$(\alpha+ \sqrt{-1}\ddbar \varphi_j)^n=
\exp(f_j)\alpha^n$$
such that $\int_X\varphi_jdV_\alpha= 0$ we have $\displaystyle 
\sup_X\vert \varphi_j\vert\leq C$.
Moreover, if $\exp(f_j)\to\exp(f_\infty)$ in $L^p$, then 
there exist a continuous function $\varphi_\infty$ such that $\varphi_j\to\varphi_\infty$ 
and which is the solution of the equation
$$(\alpha+ \sqrt{-1}\ddbar \varphi_\infty)^n=
\exp(f_\infty)\alpha^n\leqno (6)$$}

\medskip

Remark at this point that if the number $p$ above is large enough
(compared to the dimension of $X$), then the above result follows
from S.-T. Yau's original proof of the $\cC^0$ estimates (by an obvious modification), 
but the arguments provided by
S. Kolodziej seem to be more flexible, since they were adapted by
Ph. Eyssidieux, V. Guedj and A. Zeriahi in [6] (see also [11]) to get the next statement.

\medskip
\noindent{\bf Proposition {\rm [6]}} {\sl The above statement holds true 
if $\alpha\geq 0$ pointwise on $X$, 
and $\int_X\alpha^n> 0$.}
\medskip

As a 
consequence of this proposition
the family of solutions of the equation (4) admits an 
a-priori $L^\infty$ bound; this is the part of the argument where 
the integrability condition (2) on $\alpha$ is needed
(to insure the $L^p$ integrability hypothesis; remark that in the equation (*), we have
$\omega^n$ instead of $\alpha^n$ in the right hand side).

A general conclusion of the results
collected here is the next statement.

\medskip
\noindent{\bf Corollary.} {\sl There exist a constant $C> 0$ depending on $p$
and the geometry of $(X, \omega)$ such that $\vert\varphi_\vep \vert _{L^\infty}\leq C$
for each $\vep >0$. In addition, we can extract a continuous limit $\varphi$ of 
the family $(\varphi_\vep )$.
}
\medskip

\noindent {\bf Remark.} Quite recently, S. Kolodziej proved 
that the solution of the equation (6) belongs to the
H\"older space $\cC^\gamma(X)$ (for some $\gamma$ depending on $p$) 
if $\alpha$ is the inverse image of
a Kaehler metric by a generically finite map. At this moment, it is unclear whether
his methods can be applied to prove an analogous regularity result
in the hypothesis of the proposition above.

\vskip 30pt

\noindent{\twelvebf 2. End of the proof} 
\medskip

\noindent In order to achieve further regularity,
we would like to use the $\cC^2$--estimates in the theory of Monge--Amp\`ere
equations, but we cannot do it directly, because
the eigenvalues of $\alpha$ may vanish at some points of $X$.
We overcome this difficulty by using the next result of J.-P. Demailly 
and ourself (see [5]).

\medskip
\noindent{\bf Theorem {\rm [5]}} {\sl Let $(X, \omega)$ be a 
compact K\"ahler manifold, and let $\alpha$ be a semi-positive,
closed $(1,1)$ form on $X$, such that $\int_X\alpha^n> 0$.
Then there exist $\vep_0> 0$ and $\tau\in L^1(X)$ which have at worst logarithmic poles,
such that
$$\alpha+ \sqrt{-1}\ddbar \tau\geq \vep_0\omega$$ 
as currents on $X$.}

\medskip

Thus, even if the eigenvalues of the $(1,1)$--form $\alpha$ 
may be zero on an open set of $X$, we can modify it by the Hessian
of a function $\tau$, such that it dominates a small multiple of
the K\"ahler metric. Of course, now we have to deal with
the poles of $\tau$, but along the following lines we
will show that a 
careful reading of the computations
performed by Y.-T. Siu in [9] will give the result.
\vskip 5pt
\noindent Before that, remark that in general
we cannot expect the poles of $\tau$ to be divisorial
(as in the case of line bundles), so we have to
proceed to an intermediate step.
There
exist a modification (composition of blow-up maps with smooth 
centers) $\pi: \wh X\to X$ such that
$$\pi^*\alpha= \wh \omega+ [E]- \sqrt{-1}\ddbar \eta\leqno (7)$$
where $\wh \omega$ is a K\"ahler metric on $\wh X$, $E$ is an effective 
$\bQ$--divisor on $\wt X$, and $\eta$ is a quasi-psh function on $\wh X$.
Indeed, we first use a sequence of blow up maps to get 
rid of the poles of $\tau$. Thus on a 
model of $X$ the absolutely continuous part of
the inverse image of $\alpha+ \sqrt{-1}\ddbar \tau$ dominates 
a small multiple of the inverse image of a K\"ahler metric, and now
just recall the way one construct a metric on the blow-up of a manifold;
we refer to [5] for a complete description of this process.

We are going to use the equality (7) in order to study the regularity of 
$\varphi:= \lim_\vep \varphi_\vep $. We will use along the next lines the
following notation: if $f$ is a function on $X$, we denote by 
$\wh f$ the function $f\circ\pi$. 
 
On $\wh X$, the equality (5) reads as
$$\bigl(\pi^*(\alpha+ \vep\omega) + \sqrt{-1}\ddbar \wh \varphi_\vep\bigr)^n= 
(1+ \delta_\vep)
\exp \bigl(\wh \psi_{1;\vep}-\wt\psi_{2; \vep}\bigr)\Vert J(\pi)\Vert ^2\wh \omega^n\leqno (8)$$
where $\Vert J(\pi)\Vert ^2:= \omega^n/\wh \omega^n$. We denote by 
$\Phi_\vep:= \wh \varphi_\vep- \eta$, and by $\wh \omega_\vep:= \wh \omega+
\vep\pi^*\omega$; remark that the geometry of $(\wh X, \wh \omega_\vep)$
is bounded independently of $\vep$. By using the equality (7) and [8) we get
$$(\wh \omega_\vep+ \sqrt{-1}\ddbar \Phi_\vep)^n= (1+ \delta_\vep ^\prime)
\exp \bigl(\wh\psi_{1;\vep}-\wh\psi_{2;\vep}\bigr)\Vert J(\pi)\Vert ^2
\wh \omega_\vep ^n\leqno (9)$$
pointwise on $\wh X\setminus E$ (the symbols $\delta_\vep ^\prime$ are functions
which tend to zero in $\cC^\infty$ norm)

\vskip 5pt
The inequality we start with is borrowed from
Y.-T. Siu's book [9], page 99. Consider in general a compact K\"ahler manifold
$(\wh X, \wh \omega)$; we denote by $\Delta $ one half of the Laplace-Beltrami 
operator associated to the metric $\wh \omega$, and 
for each function $\Phi$ such that $\wh \omega+ \sqrt{-1}\ddbar \Phi> 0$
we denote 
by $\displaystyle \Delta _{\Phi}$
the Laplacian of the metric $\wh \omega+ \sqrt{-1}\ddbar \Phi$ on $\wh X$.

\medskip
\noindent{\bf Lemma {\rm [9]}} {\sl Let $\Phi, f\subset \cC^\infty(\wh X)$ 
be smooth functions on an
open subset $U$ of $\wh X$, such that
$$(\wh \omega + \sqrt{-1}\ddbar \Phi )^n= 
\exp (f)\wh \omega^n$$
pointwise on $U$. Then there exist a constant $C$ depending on the 
geometry of $(\wh X, \wh \omega)$ only, such that the following inequality holds true
(pointwise on the open set $U$):
$$\Delta_{\Phi}\bigl( \log(n+ \Delta\Phi
)\bigr)\geq {{1}\over {n+ \Delta \Phi}}\bigl(\Delta f- C\bigr)- 
C\sum_{j=1}^n{{1}\over{ 1+ \Phi_{, j\overline j}}}.
\leqno (10)$$}

In our case, the geometry of the family of K\"ahler manifolds
$(\wh X, \wh \omega_\vep )$ is uniformly bounded, therefore the constant 
$C$ above can be assumed to be independent of $\vep$. Also, the Hessian of 
$\wh \psi_{1;\vep }$ is bounded from below independently of 
$\vep $, thus we have $\Delta _\vep \wh \psi_{1;\vep }+ 
\log \Vert J(\pi)\Vert ^2)\geq -C$
uniformly with respect to $\vep$. Thus, the inequality (10) imply
$$\Delta_{\Phi_\vep } \log(n+ \Delta_\vep \Phi_\vep
)\geq -{{1}\over {n+ \Delta \Phi_\vep }}\bigl(\Delta_\vep 
\wh \psi_{2;\vep }+ C \bigr)- 
C\sum_{j=1}^n{{1}\over{ 1+ \Phi_{\vep, j\overline j}}}.\leqno (11)
$$
\vskip 5pt
\noindent We want next to move the term containing $\wh \psi_{2;\vep }$ inside the
$\displaystyle \Delta_{\Phi_\vep}$; for this, we need the following simple observation.

\medskip
\noindent{\bf Lemma. } {\sl There exist a constant $C> 0$ such that
$$\Delta_{\Phi_\vep }\wh \psi_{2;\vep }\geq 
{{\Delta_{\vep }\wh \psi_{2;\vep }}\over {n+ \Delta_\vep \Phi_\vep}} -
C\sum_{j=1}^n{{1}\over{ 1+ \Phi_{\vep, j\overline j}}}.\leqno (12) $$}

\medskip

\noindent {\bf Proof }of the lemma. We will prove the lemma by a local computation; 
using an appropriate coordinate system $(z^j)$ at a point $x\in \wh X$, 
the quantities under consideration are
\item {$\bullet$ }$\displaystyle \Delta_{\Phi_\vep}\wh \psi_{2; \vep }= 
\sum_j{{\wh\psi_{2;\vep, j\overline j}}\over {1+ \Phi_{\vep, j\overline j}}};$
\smallskip
\item {$\bullet$}
$\displaystyle \Delta_\vep \wh \psi_{2;\vep}= \sum_j\psi_{2;\vep, j\overline j}$.

We use again at this point the fact that the geometry of
$(\wh X, \wh \omega_\vep )$ is bounded, and thus there exist a constant $C> 0$
independent of $\vep$, such that $\sqrt{-1}\ddbar \wh \psi_{2;\vep }\geq -C\wh \omega_\vep$
on $\wh X$ (we use the statement (*) of the regularization theorem quoted above).
Therefore we have $\wh \psi_{2;\vep, j\overline j}\geq -C$, for all $j= 1,..., n$.
At the point $x$ we have

$$\eqalign {
\Delta_{\Phi_\vep }\wh \psi_{2;\vep }= & 
\sum_j{{\wh \psi_{2;\vep, j\overline j}}\over {1+ \Phi_{\vep, j\overline j}}}\cr
= & \sum_j{{\wh \psi_{2;\vep, j\overline j}+ C}\over {1+ \Phi_{\vep, j\overline j}}}- 
C\sum_j
{{1}\over{ 1+ \Phi_{\vep, j\overline j}}}\geq \cr
\geq & 
{{\Delta_{\vep }\wh \psi_{2;\vep }}\over {n+ \Delta_\vep (\Phi_\vep )}}-
C\sum_{j=1}^n{{1}\over{ 1+ \Phi_{\vep, j\overline j}}}
}$$
and thus the lemma is proved.
\vskip 5pt
\noindent
We put together the inequalities (11) and (12) and we get
$$\Delta_{\Phi_\vep}\bigl( \wh \psi_{2;\vep }+ \log(n+ \Delta_\vep (\Phi_\vep)
)\bigr)\geq -
C\sum_{j=1}^n\bigl(1+ {{1}\over{ 1+ \Phi_{\vep, j\overline j}}}\bigr)
$$
Remark also that the previous inequality and the fact that 
$$\Delta_{\Phi_\vep }(\Phi_\vep)= n- \sum_{j=1}^n{{1}\over{ 1+ \Phi_{\vep, j\overline j}}}$$
we infer 
$$\Delta_{\vep \Phi}\bigl(-2C\Phi_\vep + 
\wh \psi_{2;\vep }+ \log(n+ \Delta_\vep (\Phi_\vep)
)\bigr)\geq 
C\sum_{j=1}^n{{1}\over{ 1+ \Phi_{\vep, j\overline j}}}- C.
\leqno (13)$$

\noindent 
Recall now that $\Phi_\vep = \wh \varphi_\vep- \eta$, where the function 
$\eta $ has at most logarithmic poles along an analytic set $E$. Thus 
for each $\vep > 0$, there exist a constant $C_\vep> 0$ such that the next inequality 
holds true uniformly on $\wh X\setminus E$
$$-2C\Phi_\vep + \wh \psi_{2;\vep }+ \log(n+ \Delta_\vep \Phi_\vep)\leq C_\vep$$
(this is so because $\eta$ and $\Delta_\vep (-\eta)$ are bounded from above). 

We are now in good position to apply the maximum principle: consider 
$x_\vep  \in \wh X\setminus E$ the point where the maximum of the function
considered above
is achieved; at $x_\vep $ the relation (13) gives
$$\sum_{j=1}^n{{1}\over{ 1+ \Phi_{\vep, j\overline j}}}\leq C.\leqno (14)$$

On the other hand, the Monge-Amp\`ere equation (8) imply
$$\prod_{j=1}^n\bigl(1+ \Phi_{\vep, j\overline j}\bigr)\leq C\exp\bigl(-
\wh \psi_{2;\vep }\bigr)$$
thus for each $j=1,...,n$ we get
$$\bigl(1+ \Phi_{\vep, j\overline j}\bigr)\exp\bigl(
\wh \psi_{2;\vep }\bigr)\leq C$$
and therefore at the point $x_\vep$ we obtain 
$$\bigl(n+ \Delta_\vep \Phi_\vep\bigr)\exp\bigl(
\wh \psi_{2;\vep }\bigr)\leq C.$$

\vskip 5pt
\noindent Observe that so far, we did not used the 
uniform $L^\infty$ bound for the functions $\varphi_\vep$; 
we do it now and infer that at $x_\vep$ the next relation holds
$$\bigl(n+ \Delta_\vep \Phi_\vep\bigr)\exp\bigl(
\wh \psi_{2;\vep }- 2C\Phi_\vep\bigr)\leq C.\leqno (15)$$
Since $x_\vep$ is the maximum point of the previous function, 
we see that the inequality (15) holds true at any point of 
$X_1\setminus E$. 
\vskip 5pt

In conclusion, we have an uniform constant $C> 0$ such that 
$$n+ \Delta_\vep (\wh \varphi_\vep- \eta)\leq C\exp\bigl(
-\wh \psi_{2;\vep }+ 2C(\wh \varphi_\vep- \eta)\bigr).\leqno (14)$$
Now we have $\psi_{2;\vep}\geq \psi_2$ by the regularization theorem,
and moreover, we claim that for each $p> 0$ there exist $Y_p\subset \wh X$
such that $\displaystyle \exp\bigl(
-\wh \psi_{2}- 2C\eta\bigr)\in L^p_{loc}(\wh X\setminus Y_p)$. Indeed,
it is enough to consider the analytic set $Y_p$ where the Lelong 
numbers of the quasi-psh function $\wh \psi_{2}+ 2C\eta$ are larger than $1/p$
(for the analyticity of the level sets, see the paper of Y. T. Siu [8])
and the local integrability of the function in the complement of this set
is a consequence of a result of H. Skoda [10].

We quote now the next
regularity result (see e.g. [1])
\medskip

\noindent {\bf Theorem {\rm (Kondrakov)}} {\sl Let $B\subset \bC^n$ be the unit ball;
then the inclusion
$$L_2^p(U)\mapsto \cC^{1, \gamma}(U)$$
is compact, provided that $q(1- \gamma)> n$.}

\medskip

Thus the theorem is proved; as for the corollary, it is an easy consequence of the
previous considerations, and of the classical Schauder theory, as it is usually
applied in the context of the Monge-Amp\`ere operators (see e.g. [9], [13]) 

\vskip 10pt

\noindent {\bf Remark.} Let us consider the following geometric context:
$Y$ is a compact complex K\"ahler space (eventually singular), 
$X$ is a compact K\"ahler (smooth) manifold, $\pi: X\to Y$
is a generically finite map, $\alpha= \pi^*\omega_Y$ is the inverse image of
a K\"ahler metric on $Y$ and finally the functions $\psi_j$ have logarithmic poles.
Then it is very likely that
all the results proved here (and also in [*]) are an obvious consequence of
the original proof of the Calabi conjecture.
Indeed, a strong indication in this direction is the next observation, due to
Cascini--LaNave [2]:
the holomorphic bisectional curvature
of $\pi^*(\omega_Y)$ is bounded. On the dark side, 
the matter is not completely clear, since for example
the components of the Ricci tensor of $\pi^*(\omega_Y)$ {\it are not bounded}
(because of the additional contraction with the singular metric).

\vfill
\eject

\noindent {\twelvebf References}
\bigskip

{\eightpoint
\parindent = 1cm

\item{\bf [1]} Aubin, T\ --- {\sl Some Nonlinear Problems in Riemannian Geometry},  
Springer monographs in mathematics (1998).

\item{\bf [2]}  Cascini, P., La Nave, G.\ --- {\sl K\"ahler-Ricci Flow and the Minimal Model Program for Projective Varieties}, math.AG/0603064.

\item{\bf [3]} Demailly, J.-P.\ --- {\sl Regularization of closed positive currents and
intersection theory} J. Alg. Geom. {\bf 1}, 1992.

\item{\bf [4]} Demailly, J.-P., Koll\`ar, J.\ --- {\sl Semicontinuity of complex singularity exponents and Kähler-Einstein metrics on Fano orbifolds}, math.AG/9910118, Ann. Ec. Norm. Sup {\bf 34} (2001).

\item{\bf [5]} Demailly, J.-P., P\u aun M. \ --- {\sl Numerical characterization of the Kähler cone of a compact Kähler manifold}.  Ann. of Math. (2)  {\bf 159}  (2004),  no. 3,

\item{\bf [6]} Eyssidieux, Ph., Guedj, V., Zeriahi, A.\ --- {\sl Singular K\"ahler--Einstein metrics}, 
arXiv: math. AG/0603431.

\item{\bf [7]} Kolodziej, S. \ --- {\sl The complex Monge-Amp\`ere equation} Acta Math. {\bf 180} (1998).

\item{\bf [8]} Siu, Y.--T. \ --- {\sl Analyticity of sets associated to Lelong
numbers and the extension of closed positive currents},  Invent.\ Math.\
(1974), 53-156

\item{\bf [9]} Siu, Y.-T. \ --- {\sl Lectures on Hermitian--Einstein Metrics for Stable
Bundles and K\"ahler--Einstein metrics}, Birkh\"auser (1987).

\item{\bf [10]} Skoda, H. \ --- {\sl Sous--ensembles analytiques d'ordre fini ou infini
dans $\bC^n$}, Bull.\ Soc.\ Math.\ France, {\bf 100} (1972).

\item{\bf [11]} Tian, G., Zhang, Z. \ --- {\sl On the K\"ahler--Ricci flow of projective 
manifolds of general type}, Preprint 2005.

\item{\bf [12]} Tsuji, H. \ --- {\sl Existence and degeneration of K\"ahler--Einstein 
metrics on minimal algebraic manifolds of general type}, Math.\ Ann.\ {\bf 281} (1988).

\item{\bf [13]} Yau, S.-T.\ --- {\sl On the Ricci curvature of a 
complex K\"ahler manifold and the complex Monge--Amp\`ere equation},
Comm.\ Pure Appl.\ Math.\ {\bf 31} (1978).

\end